\documentclass[12pt]{article}                     % onecolumn (standard format)
\usepackage{amsmath}
\usepackage{bbm}
\usepackage{amsfonts,amssymb}
\usepackage[mathcal]{eucal}
\usepackage{srcltx}

\newtheorem{lemma}{Lemma}
\newtheorem{theorem}{Theorem}
\newtheorem{corollary}{Corollary}
\newtheorem{Remark}{\mbox{R e m a r k}\ \ \\ }
\newtheorem{Definition}{Definition}
\newtheorem{proposition}{Proposition}

\newcommand{\rav}{\stackrel{\triangle}{=}}
\newcommand{\valm}{\mm{V}^{\,\flat}}
\newcommand{\valp}{\mm{V}^{\,\sharp}}
\newcommand{\vplm}{\ct{V}^{\,\flat}}
\newcommand{\vplp}{\ct{V}^{\,\sharp}}

\newcommand{\ove}[3]{{#1}^{|#3\rangle}_{#2}}
\newcommand{\cp}{\mathbbm{c}}
\newcommand{\dd}[2]{\cp^{|#2\rangle}[{#1}]}
\newcommand{\rref}[1]{$(\ref{#1})$}
\newcommand{\pust}{\varnothing}
\newcommand{\mm}[1]{{\mathbb{#1}}}
\newcommand{\mA}{{\mathfrak{A}}}
\newcommand{\mB}{{\mathfrak{B}}}
\newcommand{\ct}[1]{{\mathcal{#1}}}
\newcommand{\ravref}[1]{\stackrel{(\ref{#1})}{=}}
\newcommand{\leqref}[1]{\stackrel{(\ref{#1})}{\leq}}

\newcommand{\epsi}{\varepsilon}
\newcommand{\bo}{\hfill {$\Box$}}
\begin{document}

%\setcounter{chapter}{2} % If you are doing your chapter as chapter one, %\setcounter{section}{3} % comment these %two
%lines out.

\title{\Large
On asymptotic value for dynamic games with saddle point\thanks{Supported by RFBR grant No~13-01-00304.}}
\author{Dmitry Khlopin \thanks{Krasovskii Institute  of Mathematics and Mechanics.}\\
{\it khlopin@imm.uran.ru}}
%\date{}
%\date{$\ $}

\maketitle

%\pagenumbering{arabic} %\setcounter{page}{1}%Leave this line commented out.

\begin{abstract}
The paper is concerned with two-person games with saddle point. We investigate the limits of value functions for long-time-average payoff, discounted average payoff, and the payoff that follows a probability density.

We provide two results. The first one is a uniform Tauber result for games: if the value functions for long-time-average payoff converge uniformly, then there exists the uniform limit for probability densities from a sufficiently broad set; moreover, these limits coincide. The second one is the uniform Abel result: if a uniform limit for self-similar densities exists, then the uniform limit for long-time average payoff also exists, and they coincide.

%\small\baselineskip=9pt
%The paper is concerned with two-person games with saddle point.
%    We investigate the limits of value functions
%     for long-time-average payoff,
%      discounted average payoff, and
%     the payoff that follows a probability density on $\mm{R}_{\geq 0}.$
%%\baselineskip=8pt
\end{abstract}

For bounded sequences, Hardy  proved that the convergence of their Cesaro means  is equivalent to the convergence of
their Abel means. This result has been generalized by Feller  to the case of uncontrolled deterministic dynamics in
continuous time; \cite{barton} to deterministic controlled dynamics with dependence on the initial data: if there
exists a uniform limit of values for one of the payoffs
\begin{eqnarray*}
\frac{1}{T}\int_{0}^T g\big(z(t)\big)\,dt\quad\textrm{(Cesaro mean)},\qquad\lambda\int_{0}^\infty e^{-\lambda
t}g\big(z(t)\big)\,dt\quad\textrm{(Abel mean)},
\end{eqnarray*}
then the other exists too and they are equal.

Many papers are devoted to the subject of existence of limits of such values; let us first of all note
\cite{arisconst,AG2000,BGQ2013,QG2013,QuinRenault,Vigeral2013}; see also the reviews in  \cite{barton,Ziliotto2013}.
 As proved in
 \cite{MN1981}, for a stochastic two-person game with a finite number of states and actions, optimal
 long-time~averages and optimal discounted averages share the common limit;
 for repeat games, see \cite{SIAM2012}.
For differential games, the limits exist in special cases (first of all, in the nonexpansive-like case);
\cite{AlvBarditrue,BardiGame,Cardal2010}. The general result that mirrors the result from \cite{barton}
for dynamic games (including differential
games)  was  proved in \cite{KhlopinArXiv}.

% As proved in
% \cite{MN1981}, for a stochastic two-person game with a finite number of states and actions, optimal %long-time~averages and optimal discounted averages share the common limit;
% for repeat games, see \cite{SIAM2012}.
%For differential games, the limits exist in special cases (first of all, in the nonexpansive-like case); %\cite{AlvBarditrue,BardiGame,Cardal2010}. The result that mirrors the result from \cite{barton}  for differential %games  was  announced in \cite{Khlopin1}.

In addition to the functions $\frac{1}{T}1_{[0,T]},$ $\lambda e^{-\lambda t}$, one could consider arbitrary
probability densities,  and then the asymptotic behavior of  the value functions for the following games. For
discrete time systems, \cite{Renault2013} expresses the conditions of existence of a uniform limit that is shared by
all probability distributions. In \cite{1993}, it was noted that if a probability density is nonincreasing, the
payoff can be expressed as a convex combination of the Cesaro means. Consequently, that paper proved that, for
discrete time systems, there exists a uniform limit for such distributions and it coincides with the  limit of the
long-run-average if the latter exists and is uniform. A similar approach was used in
 \cite{SIAM2012} for repeat games.

  In this paper, we show uniform Tauberian theorems for dynamic two-person games with saddle points.
  Most of our assumptions restrict the dynamics of games.
  In particular, like in \cite{barton}, we assume the closedness of strategies under concatenation.
  It is also necessary for the value function to satisfy Bellman's optimality principle, even if in a weakened,
  asymptotic sense.

 We provide two results. The first one is a uniform Tauber result for games  with saddle point:
 the uniform convergence of value functions
 for long-time average payoff implies the existence of a uniform limit of value functions
 for all
   probability densities from a sufficiently broad set (of bounded variation type); moreover, these limits
   coincide.
 The second one is the uniform Abel result for the games with saddle point:
 the existence of a uniform limit of value functions
 with self-similar
 probability densities  implies the existence of a uniform limit of value functions
 for long-time-average payoff, and these limits coincide.
% In particular, we generalize the results from \cite{barton}.

The similar results for $\frac{1}{T}1_{[0,T]},$ $\lambda e^{-\lambda t}$ were proved in \cite{KhlopinArXiv}.
Although our proofs follow \cite{KhlopinArXiv}, we use a weaker modification of Bellman's optimality principle.
This allows us to generalize the results of \cite{barton} in this article.

\section{General framework of dynamic system}
% Assume $\mm{T}\rav\mm{R}_{\geq 0}.$
  Assume the following items are given:
 \begin{itemize}
   \item   a  set  $\Omega$ of states;
   \item   a set $\mm{K}$ of feasible processes, which is a subset of mappings from $\mm{R}_{\geq 0}$ to
       $\Omega$;
   \item   a running cost $g:\Omega\mapsto [0,1];$ for each process $z\in \mm{K},$ assume the  map $t\mapsto
       g\big(z(t)\big)$ is  Borel-measurable.
 \end{itemize}

For each $\omega\in\Omega$, define $\Gamma(\omega)\rav\{z\in\mm{K}\,|\,z(0)\in\Omega\},$
    %\begin{eqnarray*}
     %\end{eqnarray*}
the set of all feasible processes $z\in\mm{K}$ starting from $\omega$.

  Let us now define the operation of  concatenation on processes.
   Let $\tau\in\mm{R}_{\geq 0},z',z''\in\mm{K}$ be such that
   $z'(\tau)=z''(0).$ Then, their concatenation $z'\diamond_\tau z''$ is defined by the following rule:
$$(z'\diamond_\tau z'') (t)\rav \left\{
 \begin{array} {rcl}        z'(t),       &\mathstrut&     t<\tau;\\
                            z''(t-\tau),
                                         &\mathstrut&     t\geq\tau.
 \end{array}            \right.$$

% {\bf Assumptions on the capabilities of players.}\
%Assume we are given a family $\mA$ of subsets of the set
%$\mm{K}$.
   Call a subset $A$ of the set $\mm{K}$ a
  {\it playable strategy} if, for every initial position $\omega\in\Omega$, we have $A\cap \Gamma(\omega)\neq\pust.$

%  Let us now define the operations on processes and playable strategies. %  For a time $\tau\in\mm{T}$ and a process
%$z\in\mm{K}$, define the function %   $z_\tau:\mm{T}\mapsto\Omega$ by the following rule: %  $$z_\tau(t)=z(t+\tau)
%\qquad\forall t\in\mm{T}.$$

  Now, for each two playable strategies $A',A''$  and a time $\tau\in\mm{R}_{\geq 0}$, define
  concatenation by
\begin{eqnarray*}
A'\diamond_\tau  A''\rav\big\{z'\diamond_\tau z''\,|\,z'\in{A}',z''\in A'',z'(\tau)=z''(0)\big\}.
\end{eqnarray*}
  Hereinafter set $A\diamond_{\tau'} A'\diamond_{\tau''} A''\rav
  \big(A\diamond_{\tau'} A'\big)\diamond_{\tau''} A''.$

Define the following axioms for some family  $\mA$ of subsets $\mm{K}$:
%every family $\mA$ of subsets $\mm{K}$:
 \begin{description}
   \item[$(p)$] $\mA$ is some family of playable strategies;
   \item[$(\omega)$] $\mA$ allows the separation of $\omega$ (at the initial time): for each mapping
   $\xi:\Omega\to \mA$
 \begin{equation*}
    %\label{300}
    \exists A\in\mA\ \forall \omega\in\Omega\qquad \xi(\omega)\cap\Gamma(\omega)=A\cap\Gamma(\omega);
 \end{equation*}
   \item[$(\diamond)$] $\mA$ is closed under concatenation $\diamond$: $\forall \tau>0,A',A''\in\mA$
   $A'\diamond_\tau A''\in\mA.$
\end{description}

\section{Formalization of lower and upper games}
 Consider a two-player game. The first player wishes  to maximize
  a bounded payoff function $c:\mm{K}\to \mm{R}$; the second player wishes to minimize it.
   The first player  also has his family $\mA$ of playable strategies specified.

  The game is conducted in the following way:
  for a given $\omega\in\Omega,$
  the first player demonstrates some   set  $A$
  in $\mA$, and
  then the second player chooses a process ${z\in A}\cap\Gamma(\omega)$;
  we have a conflict:
  $${A\in {\mA} \Uparrow}\ \ \ c(z)\ \ \ \  {\Downarrow z\in A\cap \Gamma(\omega)}.$$
  The value function of the lower game is
  \begin{eqnarray*}
 \ \ \valm[c](\omega)\rav\sup_{A\in\mA}\inf_{z\in A\cap \Gamma(\omega)}c(z)\qquad \forall\omega\in\Omega.
% \label{280}
\end{eqnarray*}
 Note that,  for every playable strategy $A$, $A\cap \Gamma(\omega)\neq\pust,$ $c$ is bounded; therefore, the value
 for the lower game is valid.

  For the upper game,
  we still have two players, the first player  maximizes
  the payoff~$c$, whereas the second player minimizes it.
  Let the second player
    also have a family $\mB$ of playable strategies.

  The upper game is conducted in the following way.
  Given $\omega\in\Omega,$
  let the second demonstrate some   set $B$
  in $\mB$;
  then, let the second player choose some process ${z\in B}\cap\Gamma(\omega)$.
  The value function of the upper game is
  \begin{eqnarray*}
 \ \ \valp[c](\omega)\rav\inf_{B\in\mB}\sup_{z\in B\cap \Gamma(\omega)}c(z)\qquad \forall\omega\in\Omega.
% \label{580}
\end{eqnarray*}
 Thanks to the boundedness of $c$ and playable strategies in $\mB$, this function is valid as well.

  Remember  that a game with the payoff function $c$
       has a {\it saddle point} if $\valm[c]=\valp[c].$
       We would also require a slightly weaker definition.
    \begin{Definition}
     Assign a bounded function $c_\lambda:\mm{K}\to\mm{R}$ to every  positive $\lambda$.
     Let us say that a game family with the payoff functions $c_\lambda$
     %probability densities $\varrho_\lambda$
       has an {\it asymptotic saddle point} if  the limit of
       $\big|\valm[c_\lambda](\omega)-\valp[c_\lambda](\omega)\big|$ as $\lambda\downarrow 0$
       exists, is equal to $0$, and
        is uniform  for $\omega\in\Omega$.
    \end{Definition}

 We hereinafter impose the following axiom on the families $\mA,\mB:$
 \begin{description}
   \item[$(s)$] $A\cap B\cap \Gamma(\omega)\neq\pust$ for any $A\in\mA,B\in\mB,\omega\in\Omega$.
 \end{description}

The utility of conditions $(s),(\omega)$ is clarified by the following lemma:
\begin{lemma}
 \label{topology}
  Suppose $\mA,\mB$ satisfy conditions $(p),(\omega)$.
  Assume $\mm{K}$ is equipped with a topology such that some bounded mapping $c:\mm{K}\to \mm{R}$ is continuous.
  In addition, let
  $cl\,A\cap cl\,B\cap \Gamma(\omega)\neq\pust$ for all $\omega\in\Omega,A\in\mA,B\in\mB.$
  Then, for all $\epsi>0$, there exist $A^\epsi\in\mA,B^\epsi\in\mB$ such that,  for all $\omega\in\Omega$,
 \begin{eqnarray}
 \valm[c](\omega)-\epsi\ <\ \inf_{z\in A^\epsi\cap \Gamma(\omega)} c(z)
 \leq\sup_{z\in B^\epsi\cap \Gamma(\omega)}c(z)<
 \valp[c](\omega)+\epsi;\nonumber\\
 \valm[c](\omega)\leq\valp[c](\omega). \label{910}
 \end{eqnarray}%

  In particular, %\rref{910} holds
   these inequalities hold if
     $\mA,\mB$ satisfy conditions $(p),(\omega),(s)$.
\end{lemma}
{\textbf{Proof}.}
 Consider arbitrary $\omega\in\Omega,\epsi>0;$  there exist maps $\xi^\epsi:\Omega\to\mA$,$\zeta^\epsi:\Omega\to\mA$
 such that, for the payoff $c$ under initial condition $z(0)=\omega,$ the strategy $\xi^\epsi(\omega)$ is
 $\epsi$-optimal in the lower game (the strategy $\zeta^\epsi(\omega)$ in the upper game).
    Then, by the property of $(\omega)$, there exist strategies $A^\epsi\in\mA$,$B^\epsi\in\mB$ for which
   $\xi^\epsi(\omega)\cap \Gamma(\omega)=A^\epsi\cap \Gamma(\omega)$,
   $\zeta^\epsi(\omega)\cap \Gamma(\omega)=B^\epsi\cap \Gamma(\omega)$ for any $\omega\in\Omega$.
 Choose a common trajectory $z^\epsi(\omega)\in cl\,A^\epsi\cap cl\,B^\epsi\cap\Gamma(\omega).$

 Since~$c$ is continuous,
 \begin{eqnarray*}
  \valm[c](\omega)-\epsi&<&\inf_{z\in cl \xi^\epsi(\omega)\cap \Gamma(\omega)} c(z)
%\inf_{z\in \xi^\epsi(\omega)\cap \Gamma(\omega)} c(z)&=&
  =\inf_{z\in A^\epsi\cap \Gamma(\omega)}c(z)\leq c\big(z^\epsi(\omega)\big)
 \end{eqnarray*}
  for all $\omega\in\Omega$.
 The converse inequality is proved for $\valp[c],B^\epsi$ in a similar way.
 Now, since the choice of $\epsi$ was arbitrary, these inequalities imply \rref{910}.

 If we know that it is always $A\cap B\cap\Gamma(\omega)\neq\pust$, then, because
 $c$ would be continuous were $\mm{K}$ equipped with the discrete topology, the proof is complete.
%\end{proof}
\bo

%\subsection{On games with asymptotic saddle point}
\section{On probability densities}

   Consider a summable Borel-measurable function $\varrho$ in $B(\mm{R}_{\geq 0},\mm{R}_{\geq 0}).$  Assign to it a
   payoff function
     $\cp[\varrho]:\mm{K}\to [0,1]$ by the following rule:
     $$\cp[\varrho](z)\rav\int_{0}^\infty \varrho (t)g\big(z(t)\big)\,dt\qquad \forall z\in\mm{K}.$$
   Note that every $\cp[\varrho]$ is bounded
    by the number $||\varrho||_{L_1(\mm{R}_{\geq 0},\mm{R}_{\geq 0})}.$
   Therefore, the following definitions are valid for all $\omega\in\Omega$:
    \begin{eqnarray*}
\vplm[\varrho](\omega)&\rav&\valm\big[\cp[\varrho]\big](\omega)=\sup_{A\in\mA}\inf_{z\in A\cap
\Gamma(\omega)}\int_{0}^\infty \varrho (t)g\big(z(t)\big)\,dt\in\mm{R},\\
\vplp[\varrho](\omega)&\rav&\valp\big[\cp[\varrho]\big](\omega)=\inf_{B\in\mB}\sup_{z\in B\cap
\Gamma(\omega)}\int_{0}^\infty \varrho (t)g\big(z(t)\big)\,dt\in\mm{R}.
    \end{eqnarray*}
  Recall that
%\begin{Definition}
  a Borel-measurable mapping $\varrho\colon\mm{R}_{\geq 0}\to\mm{R}_{\geq 0}$ is called a probability density
  function if
  $$\displaystyle\lim_{T\to\infty}\int_{0}^T\varrho(t)\,dt=1.$$
%\end{Definition}
For every density $\varrho$ and number $r\in (0,1)$ there exists the quantile $q[\varrho](r)$, i.e., the
minimum number such that
  $$\int_0^{q[\varrho](r)}\varrho(t)\,dt=r.$$

   Naturally, in this case,
       ``a game family with densities $\varrho_\lambda$
       has  an asymptotic saddle point'' means that
         the game family with the payoffs function $\cp[\varrho_\lambda]$
         has  an asymptotic saddle point.

In what follows, for every interval $[a,b)\subset\mm{R}$ and function $f:[a,b)\to\mm{R}\cup\{\infty\}$, denote by $Var_{a}^{b}[f]$ the total variation of the function $f$ in $[a,b)$. Consider the following definition:
\begin{Definition}
 We say that the family of densities $\varrho_\lambda$ ($\lambda>0$) is {\it flat at zero}
  if
      \begin{eqnarray}\label{536}\quad
  \lim_{r\to 0}\lim_{\lambda\downarrow 0} Var_{0}^{q[\varrho_\lambda](r)}[\ln \varrho_\lambda]=0.
     \end{eqnarray}
\end{Definition}
\begin{Definition}
 We say that the family of densities  $\varrho_\lambda$ ($\lambda>0$) is {\it regular }
  if, for each $r\in(0,1)$,
 the following upper limit is bounded: %THE CRUTCH дл€ борьбы с лишними строчками
      \begin{eqnarray}\label{500}\quad
  \limsup_{\lambda\downarrow0} Var_{0}^{q[\varrho_\lambda](r)}[\ln \varrho_\lambda]\in\mm{R}.
     \end{eqnarray}
\end{Definition}
%!!It is easy to verify that

Define, for every $\lambda>0$, the  probability densities $\varpi_\lambda,\pi_\lambda$ by the rule
$$\varpi_\lambda\rav\lambda 1_{[0,1/\lambda]},\qquad \pi_\lambda\rav\lambda e^{-\lambda t}.$$
For all
$r\in(0,1),\lambda>0$, we have
\begin{eqnarray*}
 Var_{0}^{q[\varpi_\lambda](r)}[\ln \varpi_\lambda]=0,\qquad  
 Var_{0}^{q[\pi_\lambda](r)}[\ln \pi_\lambda]=-\ln (1-r).
%\textrm{\ with\ }\pi_\lambda=\lambda e^{-\lambda t};
\end{eqnarray*}
Thus, the families $(\varpi_\lambda)_{\lambda>0}$ and  $(\pi_\lambda)_{\lambda>0}$ are flat at zero and regular.

For every density $\varrho$ and parameter $r\in(0,1)$, define the function
 $\ove{\varrho}{}{r}:\mm{R}_{\geq 0}\to\mm{R}_{\geq 0}$
by the rule
 \begin{eqnarray*}
% \overleftarrow{\varpi_\lambda}^r=\varpi_\lambda\cdot 1_{[0,q[\varpi_\lambda](1-r)]},\\
 \ove{\varrho}{}{r}(t)=\varrho\big(t+q[\varrho](r)\big)\qquad \forall t\geq 0.
 \end{eqnarray*}
 Let us also define the payoff $\dd{\varrho}{r}:\mm{K}\to\mm{R}$ for every $\lambda>0,r\in(0,1)$ by the following
 rule:
 for all $z\in\mm{K}$,
   $$\dd{\varrho}{r}(z)\rav
   \int_{0}^{q[\varrho](r)}   \varrho(t)g\big(z(t)\big)\,dt+
   \vplp\big[\ove{\varrho}{}{r}\big]\Big(z\big(q[\varrho](r)\big)\Big). $$

Now, for every $r\in(0,1),\lambda>0$, we have
\begin{eqnarray*}
 \ove{\varpi}{\lambda}{r}= (1-r)\varpi_{\frac{\lambda}{1-r}},\quad
 %\textrm{\ with\ }\varpi_\lambda\rav\lambda 1_{[0,1/\lambda]};\\
 \ove{\pi}{\lambda}{r}= (1-r)\pi_{\lambda}.
 %\textrm{\ with\ }\pi_\lambda\rav\lambda e^{-\lambda t}.
\end{eqnarray*}

\begin{Definition}
  The family of densities  $\varrho_\lambda$ is {\it self-similar} if, for each positive $r\in(0,1),\lambda$, there
  exists a positive $h(\lambda,r),\nu(\lambda,r)$
  such that
  \begin{equation}
  \label{right}
  \ove{\varrho}{\lambda}{r}=h(\lambda,r)\varrho_{\nu(\lambda,r)}.
  \end{equation}
\end{Definition}

Note that there are quite many such self-similar families. For example, for every Borel-measurable function
$f:\mm{R}\to\mm{R}_{>0}$, one can define its self-similar family of densities $\varrho_{\lambda}$ $(\lambda>0)$ by
the following rule:
$$\varrho_{1/T}(t)=\frac{f(T-t)1_{[0,T]}(t)}{\int_{0}^{T}f(T-\tau)\,d\tau}\qquad \forall T,t>0. $$

%!!Let us say that a family of densities has an asymptotic value function $v_*$ if ...!!

In view of the definitions of $q[\varrho_\lambda](r),\ove{\varrho}{\lambda}{r}$, condition~\rref{right} implies
 \begin{eqnarray}
 %\label{right_}
 \nonumber
1-r=\int_{0}^\infty     \ove{\varrho}{\lambda}{r}(t)dt =  \int_{0}^\infty
h(\lambda,r)\varrho_{\nu(\lambda,r)}(t)dt = h(\lambda,r),\\
 \label{right_} \ove{\varrho}{\lambda}{r}=(1-r)\varrho_{\nu(\lambda,r)}.
 \end{eqnarray}

Assume $\rho_\lambda(0)=\lambda$ for all $\lambda>0$. Then, for each $r\in(0,1)$,
$$(1-r)\nu(\lambda,r)=h(\lambda,r)\varrho_{\nu(\lambda,r)}(0)= \varrho_\lambda\big(q[\varrho_\lambda](r)\big).$$
Thus,  
  \begin{eqnarray}\label{533}
 \lim_{\lambda\downarrow 0}\nu(\lambda,r)=0\qquad\forall r\in(0,1)
 \end{eqnarray}
 if and only if
$\varrho_\lambda\big(q[\varrho_\lambda](r)\big)\to 0$
  as $\lambda\downarrow 0$ for any $r\in(0,1).$
    
 On other side, for every $r\in(0,1),$ $\lambda>0$,
$$
 \ln \varrho_\lambda\big(q[\varrho_\lambda](r)\big)\leq\ln \varrho_\lambda(0)+
 Var_{0}^{q[\varrho_\lambda](r)}[\ln \varrho_\lambda].$$
\begin{Remark}\label{751}
 If the self-similar family of $\varrho_\lambda$ is regular, and $\varrho_\lambda(0)=\lambda$ for all $\lambda>0$, then \rref{533} holds.
\end{Remark}

\section{Main result.}

\begin{theorem}
\label{1}
 Assume $\mA,\mB$ satisfy conditions $(p),(\diamond),(\omega),(s).$

  Consider a self-similar family of densities $\varrho_\lambda$ $({\lambda>0})$ such that
\begin{enumerate}
\item every function $\varrho_\lambda$ satisfies $\varrho_\lambda(0)=\lambda$;
    moreover, the family is flat at zero;
  \item  the family satisfies \rref{533}; for example, the family is regular;
  \item  for sufficiently small~$r$, the family of games with payoffs $\dd{\varrho_\lambda}{r}$ $(\lambda>0)$
  has an asymptotic saddle point;
  \item  a game family with the probability densities  $\varrho_\lambda$ $(\lambda>0)$
  has  an asymptotic saddle point, and there exists a
  common limit $V_*$ for  $\vplm[\varrho_\lambda]$ and $\vplp[\varrho_\lambda]$ that is uniform in
$\omega\in\Omega$.
\end{enumerate}

 Then, the game family with probability densities $\varpi_\mu=\mu 1_{[0,1/\mu]}$ ($\mu>0$) has
  an asymptotic saddle point,  and the common limit of values $\vplp[\varpi_\mu]$ and $\vplp[\varpi_\mu]$ also
  exists, is uniform in
$\omega\in\Omega$, and coincides with~$V_*.$
\end{theorem}

\begin{theorem}
\label{2}
  Assume $\mA,\mB$ satisfy axioms $(p),(\diamond),(\omega),(s).$

  Let a regular family of densities $\varsigma_\mu$ $(\mu>0)$ satisfy
%  In what follows, we are only concerned with the families of densities~$\varrho_\lambda$ for which
\begin{equation}
\label{520}
   \lim_{\mu\to 0} \int_{0}^{T}\varsigma_\mu(t)\,dt=0 \qquad\forall T>0.
\end{equation}

  For the family with probability densities $\varpi_\lambda=\lambda 1_{[0,1/\lambda]}$ $(\lambda>0)$ assume:
\begin{enumerate}
  \item  under sufficiently small~$r$, the game family with the payoffs $\dd{\varpi_\lambda}{r}$ $(\lambda>0)$
      has an asymptotic saddle point;
  \item  a game family with probability densities $\varpi_\lambda$
  has  an asymptotic saddle point, and there exists a common limit $V_*$ for  $\vplm[\varpi_\lambda]$ and
  $\vplp[\varpi_\lambda]$ that is uniform in
$\omega\in\Omega$.
\end{enumerate}

 Then,
a game family with probability densities $\varsigma_\mu$
  has  an asymptotic saddle point,  the common limit of payoffs $\vplp[\varsigma_\mu]$ and $\vplp[\varsigma_\mu]$
  also exists, is uniform in
$\omega\in\Omega$, and coincides with~$V_*.$
\end{theorem}

%\section{Uniform Tauber-Abel theorem if saddle point always exists}

\begin{corollary}
\label{4}
 Assume $\mA,\mB$ satisfy conditions $(p),(\omega),(\diamond),(s).$

  Let a game family with payoffs $\dd{\varrho_\lambda}{r}$ $({\lambda>0})$ have an asymptotic saddle point for every
  family of densities $\varsigma_\lambda$ and all
$r\in(0,1)$.

 Then, the following conditions on the function $V_*:\Omega\to[0,1]$ are equivalent:
         \begin{itemize}
  \item\  for a certain regular self-similar family of densities $\varrho_\lambda$ with
      $\varrho_\lambda(0)=\lambda$ that is flat at zero and satisfies \rref{520},
 there exists the common limit
\begin{equation}\label{1113}
  \lim_{\lambda\downarrow 0} \vplp[\varrho_\lambda](\omega)=\lim_{\lambda\downarrow 0}
  \vplm[\varrho_\lambda](\omega)=V_*(\omega);
  \end{equation} that is uniform for $\omega\in\Omega;$
  \item \  for a  family of densities $\varrho_\lambda=\lambda 1_{[0,1/\lambda]}$, there exists the common limit of
      \rref{1113}
 that is uniform for $\omega\in\Omega;$
 \item\  for a family of densities $\varrho_\lambda=\lambda e^{-\lambda t} $, there exists the common limit of
     \rref{1113}
 that is uniform for $\omega\in\Omega;$
\item\  for every regular family of densities $\varrho_\lambda$ $(\lambda>0)$ with \rref{520},
 there exists the common limit of \rref{1113}
 that is uniform for $\omega\in\Omega.$
         \end{itemize}
\end{corollary}

Note that the existence of an asymptotic saddle point for a family of games with the payoffs
$\dd{\varrho_\lambda}{r}$ $(\lambda>0)$ can be demonstrated in different ways. For example, under conditions
$(p),(\omega),(\diamond),(s)$, it is sufficient to prove that for every $T>0$, any $A\in\mA,B\in\mB$ can be expressed
in the form $A=A\diamond_\tau A'$, $B=B\diamond_\tau B'$
  for certain $A'\in\mA,B'\in\mB$ (see \cite{KhlopinArXiv}).

  In some cases, not even that is necessary.

{\bf Case of one-person games.}

Let us consider the example in the form of
 deterministic dynamic programming problem in continuous time.
 As in \cite{barton}, assume the sets
$\Omega,\mm{K}$ to be given; moreover, assume  $\mm{K}$ to be closed under concatenation.

 Consider an arbitrary bounded
  payoff function $c:\mm{K}\to\mm{R}.$
 Set $\mA\rav\{\mm{K}\}$.
 Now, %for every $\omega\in\Omega$,
 \begin{eqnarray*}
 \valm[c](\omega)=\sup_{A\in\{\mm{K}\}}\inf_{z\in A\cap\Gamma(\omega)}c(z)=
 \inf_{z\in\Gamma(\omega)}c(z)\quad\forall\omega\in\Omega.
 %\int_{0}^\infty \varrho_\lambda(t)g(z(t))\,dt\\
 %&=&\inf_{z\in\Gamma(\omega)}\int_{0}^\infty \varrho_\lambda(t)g(z(t))\,dt.
 \end{eqnarray*}
 Consider all the possible selectors $\zeta$ of the multivalued mapping $\Omega\ni\omega\mapsto\Gamma(\omega)\subset\mm{K};$
 let $\mB$ be the set of all possible images $\zeta(\Omega)$ of these mappings.
 Since $\mm{K}$ is closed under concatenation, condition $(\diamond)$ holds for~$\mA$ and~$\mB$ defined in such way.
 Conditions $(p),(\omega),(s)$ can be verified directly.
  Now, Lemma~\ref{topology} implies that
 \begin{eqnarray*}
 \valp[c](\omega)=\sup_{B\in\mB}\inf_{z\in B\cap\Gamma(\omega)}c(z)=
 %\int_{0}^\infty \varrho_\lambda(t)g(z(t))\,dt\\
 %&=&
 \inf_{z\in\Gamma(\omega)}c(z) \quad  \forall \omega\in\Omega.
 %\int_{0}^\infty \varrho_\lambda(t)g(z(t))\,dt.
 \end{eqnarray*}
 Thus, a game with arbitrary bounded payoff~$c$ has a saddle point. Finally,
\begin{Remark}
\label{5}
 For a one-person game under conditions of \cite{barton}, result of Corollary~\ref{4} holds.
\end{Remark}
In particular, the condition of uniformity for the limit of \rref{1113} is in the general case indispensable for the
results above; for the counterexample, refer to \cite{barton}.

\section{Proof of Theorems~\ref{1},\ref{2}.}
\subsection{Auxiliary statements.}

\begin{Remark}
\label{825}
  Under conditions of Theorem~\ref{1}, in view of \rref{right},\rref{right_}, for all
  $z\in\mm{K}$,
   $\dd{\varrho_\lambda}{r}(z)$ coincides with
   $$
   \int_{0}^{q[\varrho_\lambda](r)}    \varrho_\lambda(t)g\big(z(t)\big)\,dt+
   (1-r)\vplp\big[\varrho_{\nu(\lambda,r)}\big]\big(z(q[\varrho_\lambda](r))\big);$$
   moreover, by  \rref{533},
   the family of games with payoffs
   $$
   \int_{0}^{q[\varrho_\lambda](r)}    \varrho_\lambda(t)g\big(z(t)\big)\,dt+
   (1-r)\vplm\big[\varrho_{\nu(\lambda,r)}\big]\big(z(q[\varrho_\lambda](r))\big)$$
   has an asymptotic saddle point for sufficiently small~positive $r.$
\end{Remark}

\begin{lemma}
 \label{equa}
  Assume all conditions of Theorem~\ref{1} hold. Then, for all sufficiently small  $r\in(0,1)$, there exist the
  limits
    $$\lim_{\lambda\downarrow 0} \valm[\dd{\varrho_\lambda}{r}](\omega)=\lim_{\lambda\downarrow 0}
    \valp[\dd{\varrho_\lambda}{r}](\omega)=V^*(\omega),$$
    and these limits are uniform for $\omega\in\Omega.$
\end{lemma}
{\textbf{Proof}.}
  By the condition, the family of games with the payoffs $\dd{\varrho_\lambda}{r}$ has an asymptotic saddle point for
  all sufficiently small $r\in(0,1)$. Fix one such $r\in(0,1).$

   Let us show that, for all positive $\epsi<1$, $\valp[\dd{\varrho_\lambda}{r}]\geq V_*-\epsi$
    for sufficiently small $\lambda>0$.
   Suppose it is not; then, there exists a positive $\delta>0$ such that, for every natural $n$,  there exist a
   positive $\lambda_n<1/n$
    and $\omega_n\in\Omega$ such that
  $\valp[\dd{\varrho_{\lambda_n}}{r}](\omega_n)<V_*(\omega_n)-8\delta.$
  Fix that $\delta<1/8$.

 But, there exists a natural $n$ such that
  \begin{eqnarray*}
 \vplp[\varrho_{\nu(\lambda_n,r)}]<V_*+\delta,\qquad
   \vplm[\varrho_{\lambda_n}]>V_*-\delta.
 \end{eqnarray*}
Fix such $n,\lambda_n,\omega_n$.

  Now, there exists $A\in\mA$ such that %, for all $z\in A$, %
  we obtain
\begin{eqnarray}
\label{64A}
  V_*(z(0))-2\delta
  <\vplm\big[\varrho_{\lambda_n}\big]\big(z(0)\big)-\delta<
   \int_{0}^{\infty}  \varrho_{\lambda_n}(t)g\big(z(t)\big)\,dt\qquad  \forall z\in A;
\end{eqnarray}
   there exists $B''\in\mB$ such that, for all $z\in B''$, we obtain
\begin{eqnarray}
\label{64B2}
  \vplp\big[\varrho_{\nu(\lambda_n,r)}\big]\big(z(0)\big)+\delta
      >
   \int_{0}^{\infty}  \varrho_{\nu(\lambda_n,r)}(t)g\big(z(t)\big)\,dt
       \ravref{right_}
   \frac{1}{1-r}\int_{0}^{\infty}  \ove{\varrho}{{\lambda_n}}{r}(t)g\big(z(t)\big)\,dt.
\end{eqnarray}
 By definition of $\delta$, there exists $B'\in\mB$ such that, for all $z\in B'\cap\Gamma(\omega_n)$,
  we have
\begin{eqnarray}
\label{64B1}
  \dd{\varrho_{\lambda_n}}{r}(z)<\valp[\dd{\varrho_{\lambda_n}}{r}](\omega_n)+\delta<V_*(\omega_n)-7\delta.
\end{eqnarray}
   Set $T\rav q[\varrho_{\lambda_n}](r);$
 since the left-hand side of this inequality depends only on $z|_{[0,T]},$
 the strategy $B'$ can be substituted here with strategy $B'\diamond_{T}B''.$

 Then, for any  $z\in A\cap(B'\diamond_{T}B'')\cap\Gamma(\omega_n),$
\begin{eqnarray*}
   V_*(\omega_n)-2\delta&\leqref{64A}& \int_{0}^{\infty}  \varrho_{\lambda_n}(t)g\big(z(t)\big)\,dt\\
   &=&\int_{0}^{T}  \varrho_{\lambda_n}(t)g\big(z(t)\big)\,dt+
   \int_{T}^{\infty}  \ove{\varrho}{{\lambda_n}}{r}(t-T)g\big(z(t)\big)\,dt\\
   &\leqref{64B2}&\dd{\varrho_{\lambda_n}}{r}(z)+\delta
   \leqref{64B1}V_*(\omega_n)-6\delta.
\end{eqnarray*}

  The obtained contradiction proves that $\valp[\dd{\varrho_\lambda}{r}]\geq V_*-\epsi$ for any $\epsi>0$ and
  sufficiently small $\lambda>0$.
    The inequality $\valm[\dd{\varrho_\lambda}{r}]\leq V_*+\epsi$, in view of Remark \ref{825}, is proved similarly.
  It remains to note that, by the condition,
  $\valp[\dd{\varrho_\lambda}{r}](\omega)-\valm[\dd{\varrho_\lambda}{r}](\omega)\to 0$ as $\lambda\downarrow 0$, and this
  limit is uniform in $\omega\in\Omega.$
\bo
%\end{proof}

\subsection{Proof of Theorem~\ref{1}.}

In~Section \ref{233}, we proved the following proposition:
\begin{proposition}
\label{tochi}
  Let all conditions of Theorem~\ref{1} hold.
   Then, for all positive $\epsi$, there exists $\bar\mu>0$ such that, for all positive $\mu<\bar\mu$ for all
   $\omega\in\Omega$, we have
      $\vplm[\varpi_\mu](\omega)>V_*(\omega)-6\epsi.$
\end{proposition}

 Define a function $g^-\rav 1-g.$
 It is easy to see that, for a probability density $\varrho$, relation $g+g^-\equiv 1$ implies that, for all
 $\omega\in\Omega$,
 \begin{eqnarray*}
 1-\vplp[\varrho](\omega)=\inf_{B\in\mB}\sup_{z\in B} \int_{0}^\infty \varrho(t) g^-\big(z(t)\big)\,dt.
\end{eqnarray*}
 Conditions of Proposition~\ref{tochi} hold for the lower game with the probability density $\varpi_\mu,$  running
 cost $g^-$,   the limit value $1-V_*$, and swapped capabilities of players $\mB,\mA$.
 Then, for all positive $\epsi$ for sufficiently small $\mu>0$, we have
      $1-\vplp[\varpi_\mu]>1-V_*-6\epsi.$
 Now, in view of Lemma \ref{topology}, we have
      $\vplm[\varpi_\mu]\leq\vplp[\varpi_\mu]<V_*+6\epsi.$
 As proved above, for sufficiently small positive~$\mu$, we have
 $V_*-6\epsi<\vplm[\varpi_\mu].$
 Since the positive $\epsi$ was chosen arbitrarily, $\vplm[\varpi_\mu]$ and $\vplp[\varpi_\mu]$ have the common limit
 (as $\mu\to 0$); this limit equals
$V_*$ and is uniform in $\omega\in\Omega.$
\bo

%\section{Proof of Proposition \ref{tochi}.} %\label{373}
\subsection{Proof of Theorem~\ref{2}.}

Note that the family of~$\varpi_\lambda$ satisfies all conditions of Theorem~\ref{1}. Then, Lemma~\ref{equa} and
Remark~\ref{825} hold for it.

In Section \ref{373}, we proved the following proposition:
\begin{proposition}
\label{chito}
  Let all conditions of Theorem~\ref{2} hold.
   Then, for all positive $\epsi$, there exists $\bar\mu>0$ such that, for all positive $\mu<\bar\mu$ for all
   $\omega\in\Omega$, we have
      $\vplm[\varsigma_\mu](\omega)>V_*(\omega)-6\epsi.$
\end{proposition}

Now, a literal repetition of the proof of Theorem~\ref{1} with
 $\varpi_\mu$ replaced with $\varsigma_\mu$ and $\varrho_\lambda$ replaced with $\varpi_\lambda$ provides what was
 required.
\bo

\section{Proof of Proposition \ref{tochi}.}
\label{233}

   According to the condition of Theorem~\ref{1} and the result of Lemma~\ref{equa},
   for all positive $\varkappa$, for all $r\in(0,1)$, there exists a positive
   $\hat\lambda(\varkappa,r)$ such that, for all positive
   $\lambda\leq\hat\lambda(\varkappa,r)$,
\begin{eqnarray*}
  V_*-\valm\big[\dd{\varrho_{\lambda}}{r}\big]<\varkappa,\qquad
   \vplm[\varrho_{\lambda}]-V_*<\varkappa.
 \end{eqnarray*}

  Let us show that, for all positive $\epsi<1/4$, there exists $\bar\mu>0$ such that, for all positive $\mu<\bar\mu$,
  we have
      $\vplm[\varpi_\mu]>V_*-4\epsi.$
  Fix a positive $\epsi<1/4.$

By \rref{536}, for some $\check{r}\in(0,1)$, every positive $r<\check{r}$
has a positive $\check\lambda(r)$ such that, for a positive $\lambda<\check\lambda(r)$,
 $$Var_{0}^{q[\varrho_\lambda](r)}\big[\ln \varrho_{\lambda}\big]<\frac{\epsi}{2};$$
  now, for all nonnegative $t\leq q[\varrho_\lambda](r)$,
   $$\varrho_{\lambda}(t)\leq e^{\epsi/2}\varrho_{\lambda}(0)=e^{\epsi/2}\lambda<(1+\epsi)\lambda,$$
    \begin{equation}
   \label{250}
   \int_{0}^{\tau} |\varrho_{\lambda}(t)-\lambda|\, dt\leq \epsi\lambda\tau
   \end{equation}
   as soon as $\lambda<\check\lambda(r)$,$\tau\leq q[\varrho_\lambda](r), r<\check{r}$ for some positive
   $\check{r}.$ Fix a positive $\check{r}.$

% Remember that $\varpi_\mu=\mu 1_{[1,1/\mu]}.$

%$\ln k<2k\epsi,\epsi<1,$
     There exists a natural $k>2/\epsi>4$ such that $ \frac{\ln k}{k}<\check{r}$. Fix~$k$.
Set
$$p\rav k^{-1/k},\quad \varkappa\rav\frac{1}{3k^2}.$$
  By the condition,
 \begin{eqnarray}
  \label{275}
  \ \ p^{k}=\frac{1}{k},\qquad 1-p<-\ln p=\frac{\ln k}{k}<\check{r}.
 \end{eqnarray}

   From  \rref{533}, we know that there exists
a positive $\bar\lambda(k)$ such that, for positive $\lambda<\bar\lambda(k)$, we have
      $\lambda,\nu(\lambda,1-p)<\hat\lambda(\varkappa,1-p)$.
%            \vplm[\varrho_{\lambda}]-\varkappa<\valm[\dd{\varrho_{\lambda}}{1-p}].$$
   Assume $\bar{\mu}=\min\big(\bar\lambda(k),\check\lambda(1-p)\big)/k.$
   Fix $\mu<\bar\mu.$

Set $\tau_0\rav 0$;  for each $m=1,2\dots,k$, set
   $$
   \lambda_m\rav\frac{\mu p^{m-1}}{1+\epsi},\
   t_m\rav q[\varrho_{\lambda_m}](1-p),\ \tau_m\rav\tau_{m-1}+t_m.$$
  Now, in accordance with \rref{275}, we have $\lambda_m<k\mu\leq k\bar{\mu}\leq\bar\lambda(k)$; hence,
  $\lambda_m,\nu(\lambda_m,1-p)<\hat\lambda(\varkappa,1-p)$, i.e.,
 \begin{eqnarray}
 %   \lambda_m\leqref{255} \mu\delta^{-2k}
 %   \leq
 %   \mu k<\breve\lambda(k),\nonumber\\
 %   r_m\in[\alpha_k,\beta_k],\nonumber\\
   V_*-\valm[\dd{\varrho_{\lambda}}{1-p}]<\varkappa,\label{317}\\
   \vplm[\varrho_{\nu(\lambda_m,1-p)}]-V_*<\varkappa.\label{316}
 \end{eqnarray}
  Moreover, by \rref{275}, $1-p<\check{r}$, $\lambda_m<k\bar{\mu}\leq\check\lambda(1-p)$ and
$t_m=q[\lambda_m](1-p)$, we have \rref{250}, i.e.,
     \begin{equation}
  \label{251}
   \int_{0}^{t_m} |\varrho_{\lambda_m}(t)-\lambda_m|\, dt\leq \epsi\lambda_m t_m.
   \end{equation}
  Now, for all $z\in\mm{K}$,
     \begin{eqnarray}
 \nonumber
 \int_{0}^{t_m} \varrho_{\lambda_m}(t)g\big(z(t)\big)\,dt&\leqref{251}&
 (1+\epsi)\int_{0}^{t_m}\lambda_m g\big(z(t)\big)\,dt\\
 \label{253}
 &=&\mu p^{m-1}\int_{0}^{t_m}g\big(z(t)\big)\,dt;
   \end{eqnarray}
   moreover, using $\epsi<1/4,$ $t_m=q[\varrho_{\lambda_m}](1-p)$ we obtain
   \begin{eqnarray*}
 %  \nonumber
     t_m\mu p^{m-1}&=&t_m\lambda_m(1+\epsi)\\
     &\leqref{251}&
   \frac{1+\epsi}{1-\epsi}\int_{0}^{t_m}\varrho_{\lambda_m}(t)\,dt
  =\frac{(1-p)(1+\epsi)}{1-\epsi}\\
  &<& (1-p)(1+3\epsi).
 %  \nonumber
   \end{eqnarray*}
 Thus, we have $\tau_k\mu\leq(1-p)(1+3\epsi)\big(1+\dots+p^{k-1}\big)< 1+3\epsi,$ and
  \begin{eqnarray}
    \label{252}
     %implies $
     \mu 1_{[0,\tau_k]}\leq\varpi_{\mu}+\mu 1_{(1/\mu,(1+3\epsi)/\mu]}.
  \end{eqnarray}

    By \rref{317}, $V_*-\varkappa$ does not exceed the
    lower value $\valm[\dd{\varrho_{\lambda_m}}{1-p}]$ of the game with the payoff $\dd{\varrho_{\lambda_m}}{1-p},$
    i.e., with the payoff
   $$%\sup_{A\in\mA}\inf_{z\in A\cap \Gamma(\omega)}\Big[
   \int_{0}^{t_m}\varrho_{\lambda_m}(t)g\big(z(t)\big)\,dt+p
   \vplm[\varrho_{\nu(\lambda_m,1-p)}]\big(z(t_m)\big).
   %\Big].\
   $$
   In view of \rref{316}  this, in its own turn, does not exceed the
    lower value of the game with the payoff
   $$%\sup_{A\in\mA}\inf_{z\in A\cap \Gamma(\omega)}\Big[
   \int_{0}^{t_m}\varrho_{\lambda_m}(t)g\big(z(t)\big)\,dt+p
   V_*\big(z(t_m)\big)+\varkappa.
   %\Big].\
   $$
   By Lemma~\ref{topology}, this game has a $\varkappa$-optimal strategy $A^m\in\mA$; then, for all $z\in A^m$,
   $V_*\big(z(0)\big)$ does not exceed
 \begin{eqnarray*}
 \int_{0}^{t_m}\varrho_{\lambda_m}(t)g\big(z(t)\big)\,dt+p
   V_*\big(z(t_m)\big)+3\varkappa
   \leqref{253}
   p^{m-1}\int_{0}^{t_m}\mu g\big(z(t)\big)\,dt
    +p
   V_*\big(z(t_m)\big)+3\varkappa.
 \end{eqnarray*}
   Since the right-hand sides of these inequalities depend only on $z|_{[0,t_m]},$
 the strategy $A^m$ can be substituted here with any strategy
  that may be represented in the form
  $A^m\diamond_{t_m}A''.$
   Now, for all  $m=1,\dots,k,$ $A',A''\in\mA,z\in A'\diamond_{\tau_{m-1}}A^m\diamond_{\tau_{m}}A''$,
   we obtain %that $V_*\big(z(\tau_{m-1})\big)$ does not exceed
\begin{equation}
\label{369}
V_*\big(z(\tau_{m-1})\big)\leq
 p^{m-1}  \int_{\tau_{m-1}}^{\tau_{m}}    \mu g\big(z(t)\big)\,dt
    +p
   V_*\big(z(\tau_{m})\big)+3\varkappa.
\end{equation}

   Set $A^*\rav A^1\diamond_{\tau_1}A^2\diamond_{\tau_2}\dots \diamond_{\tau_{k-1}}A^k.$
    For every $z\in A^*$, \rref{369} holds for all $m=1,\dots,k.$
      Now, taking at first $m=1$, and subsequently using
\rref{369} for $m=2,3\dots,k$, we have
\begin{eqnarray*}
  V_*\big(z(0)\big)&\leqref{369}&\int_{0}^{\tau_1}\mu g\big(z(t)\big)\,dt
    +p
   V_*\big(z(\tau_1)\big)+3\varkappa\\
&\leqref{369}&\int_{0}^{\tau_2}\mu g\big(z(t)\big)\,dt
    +p^2
   V_*\big(z(\tau_2)\big)+6\varkappa\\
   \leqref{369}\dots
   &\leqref{369}&
   \int_{0}^{\tau_k}\mu g\big(z(t)\big)\,dt
    +p^k
   V_*\big(z(\tau_k)\big)+3k\varkappa\\
   &\leqref{275}&
   \int_{0}^{\tau_k}\mu g\big(z(t)\big)\,dt
    +\frac{2}{k}.
 \end{eqnarray*}
 But then, $V_*\leq \mu \vplm[1_{[0,\tau_k]}]+\frac{2}{k}$ as well.
 Now, by \rref{252},
  we have
 \begin{eqnarray*}
\mu \vplm[1_{[0,\tau_k]}]&\leq& \vplm[\varpi_{\mu}]+3\epsi,\\ V_*&\leq& \vplm[\varpi_{\mu}]+3\epsi+\frac{2}{k}.
 \end{eqnarray*}
 Since by the choice of~$k$ we have $\epsi k>2$, we obtain $V_*\leq \vplm[\varpi_{\mu}]+4\epsi$
 for all $\mu<\bar\mu,$ which was to be proved.
\bo

\section{Proof of Proposition \ref{chito}.}
\label{373} % Suppose $\varsigma_\lambda=\lambda 1_{[1,1/\lambda]}.$
  Let us show that, for every positive $\epsi<1$, there exist $\bar\mu>0$ such that, for all positive $\mu<\bar\mu$
  for all $\omega\in\Omega$, we have
      $$\vplm[\varsigma_\mu](\omega)>V_*(\omega)-6\epsi.$$
  Fix a positive $\epsi<1.$

   For every positive $\mu$, define $T(\mu)\rav q[\varsigma_{\mu}](1-\epsi)$ by $\int_{0}^{T(\mu)}
   \varsigma_{\mu}(t)\,dt=1-\epsi.$
   By \rref{500}, there also exist positive $M,\breve{\mu}$ such that
   $$Var_{0}^{T(\mu)}[\ln \varsigma_{\mu}]\,dt<M$$
   for $\mu<\breve{\mu}$.
   Fix $M,\breve{\mu}.$

   Since, for all $t\in\big[0,T(\mu)\big]$, we have
   $\varsigma_\mu(0)\leq e^M\varsigma_\mu(t),$ integrating over this interval, we obtain
   $\varsigma_\mu(0)<e^M/T(\mu)$. Then,  $\varsigma_\mu(t)<e^{2M}/T(\mu).$
   Thanks to \rref{520}, we have $T(\mu)\to\infty$ as $\mu\downarrow 0$.
   Thus, for every positive $\lambda$, there exists a positive $\hat{\mu}(\lambda)<\breve{\mu},$
   such that, for positive $\mu<\hat{\mu}(\lambda)$, we obtain $\varsigma_\mu(t)<\lambda$
   for all $t\in\big[0,T(\mu)\big].$

      There exists a natural $k>M$ such that $e^{M/k}<1+\epsi,$
       $k\epsi>-\ln \epsi.$ Then,  $k(1-\epsi^{1/k^2})<\epsi.$
       Fix such $k$.

  Set
   \begin{equation}
   \label{M}
  \varkappa\rav\frac{\epsi}{3k^2},\qquad p\rav\epsi^{1/k^2},\qquad \delta\rav 1-p<\frac{\epsi}{k}.
   \end{equation}
 Then,   $\delta+p\delta+\dots+p^{k^2-1}\delta=1-\epsi,$
  $\delta+p\delta+p^2\delta+\dots=1.$

     According to the condition of Theorem~\ref{2} and the result of Lemma~\ref{equa},
 for a certain positive $\hat{\lambda}(k)$ for all positive
     $\lambda<\epsi\hat\lambda(k)$, we have
   \begin{equation}
   \label{1053}
   \vplm[\varpi_\lambda]-V_*<\varkappa,\qquad V_*-\valm\big[\dd{\varpi_{\lambda}}{1-p}\big]<\varkappa.
   \end{equation}

   Then, by definition of $\hat{\mu}$, for some positive $\bar{\mu}<\hat{\mu}\big(\epsi\hat{\lambda}(k)\big)$, we also have
   $\varsigma_{\mu}(t)<\epsi\hat\lambda(k)$ for all $t\in\big[0,T(\mu)\big]$, $\mu<\bar{\mu}$.
      Fix such $\mu>0$.

  Note that $T(\mu)=q[\varsigma_{\mu}](1-p^{k^2}).$
  Let us partition the interval $\big[0,T(\mu)\big]=\big[0,q[\varsigma_{\mu}](1-p^{k^2})\big]$ into $k^2$ nonempty intervals
  $[\tau_m,\tau_{m+1})$ ($m=1,\dots,k^2$);
  to this end, define
 \begin{eqnarray*}
   \tau_0=0,\quad\tau_{m}=q[\varsigma_{\mu}](1-p^m),\quad
   \lambda_{m}\rav \frac{\int_{\tau_{m-1}}^{\tau_{m}}
   \varsigma_{\mu}(t)\,dt}{\tau_{m}-\tau_{m-1}}=\frac{p^{m-1}-p^{m}}{\tau_{m}-\tau_{m-1}}.
   %=\frac{p^{m-1}\delta}{\tau_{m}-\tau_{m-1}}.
 \end{eqnarray*}
   Note that the choice of~$\mu$ implies that
   $\lambda_{m}<\epsi\hat\lambda(k)$ for all
   $m=1,\dots,k^2$.
  Then, for all
   $m=1,\dots,k^2$,
   $${\lambda_m p^{1-m}}<{\lambda_m p^{-m}}<{\lambda_m p^{-k^2}}=
  \lambda_m/\epsi<\hat\lambda(k),$$
  which implies the validity of \rref{1053} for $\lambda_m p^{1-m}$ and $\lambda_m p^{-m}.$

  Call an interval $[\tau',\tau'')\subset \mm{R}$ correct if
 $Var_{\tau'}^{\tau''}[\ln \varsigma_{\mu}]\,dt<M/k.$

  Let us define a scalar function $\varsigma'$ on $\mm{R}_{\geq 0}$  by the following rules: $\varsigma'(t)=0$ for
  all $t\geq T(\mu)$; $\varsigma'(t)=\varsigma_\mu(t)$ for all $t\in[\tau_{m-1},\tau_{m})$ if this
  interval is correct, and
     $\varsigma'(t)=\lambda_m$ for $t\in[\tau_{m-1},\tau_{m})$ if this interval is incorrect.

  Note that in this case there are at least $k^2-k$ correct intervals of $k^2$.
  Then, the integrals of the functions $\varsigma_\mu$ and $\varsigma'$  over the incorrect intervals do not exceed
  $k\delta<\epsi$ by \rref{M}.
  Now, for all $z\in\mm{K}$, we have
  $$\int_{0}^\infty \varsigma'(t)g\big(z(t)\big)\,dt-
  \int_{0}^\infty \varsigma_\mu(t)g\big(z(t)\big)\,dt<2\epsi.$$
  In particular,
\begin{equation}
\label{1078}
  \vplm[\varsigma_\mu]>\vplm[\varsigma']-2\epsi.
\end{equation}

  Note that $\varsigma'$ was introduced such that, for all $m=1,\dots,k^2,$
  %$\tau_{m}=q[\varsigma'](1-p^m)]$,
   we have
\begin{eqnarray}
    \lambda_{m}=\frac{\int_{\tau_{m-1}}^{\tau_{m}}
    \varsigma'(t)\,dt}{\tau_{m}-\tau_{m-1}}=\frac{\int_{\tau_{m-1}}^{\tau_{m}}
    \varsigma_\mu(t)\,dt}{\tau_{m}-\tau_{m-1}}=\frac{p^{m-1}-p^m}{\tau_{m}-\tau_{m-1}},\nonumber\\
   \lambda_{m}\leq
   \sup_{t\in[\tau_{m-1},\tau_{m})}
   \varsigma'(t)\leq e^{M/k}\inf_{t\in[\tau_{m-1},\tau_{m})}\varsigma'(t),\nonumber\\
\label{509}
   \int_{\tau_{m-1}}^{\tau_{m}}     \lambda_{m}g\big(z(t)\big)dt\leq  e^{M/k}
     \int_{\tau_{m-1}}^{\tau_{m}}    \varsigma'(t)g\big(z(t)\big)dt.
\end{eqnarray}

 Remember that $\ove{\varpi}{{\nu}}{1-p}=p\varpi_{\nu/p}$
 for all positive $\nu.$
 Now, for all $z\in\mm{K},$  the density  $\varpi_{\lambda_m p^{1-m}}$ satisfies
\begin{eqnarray}
\nonumber  \ove{\varpi}{{\lambda_m p^{1-m}}}{1-p}=p\varpi_{\lambda_m p^{-m}},\\ \nonumber  q[\varpi_{\lambda_m
p^{1-m}}](1-p)=%\frac{1-p}{\lambda_m p^{1-m}}=
  \frac{1-p}{\lambda_m p^{1-m}}=\tau_m-\tau_{m-1},\\
\dd{\varpi_{\lambda_m  p^{1-m}}}{1-p}(z)=
  \int_{0}^{\tau_m-\tau_{m-1}}    \lambda_m p^{1-m} g\big(z(t)\big)\,dt+
\label{1027}   p\vplm\big[\varpi_{\lambda_m p^{-m}}\big]\big(z(\tau_m-\tau_{m-1})\big).
\end{eqnarray}

   Consider any of $m=1,\dots,k^2.$
    By \rref{1053},  $p^{m-1}V_*-\varkappa$ does not exceed the
    lower value of the game with the payoff $p^{m-1}\dd{\varpi_{\lambda_m  p^{1-m}}}{1-p},$ i.e., (see \rref{1027}),
    the payoff
   $$%\sup_{A\in\mA}\inf_{z\in A\cap \Gamma(\omega)}\Big[
  \int_{0}^{\tau_m-\tau_{m-1}}  \lambda_m  g\big(z(t)\big)\,dt+
   p^m\vplm\big[\varpi_{\lambda_m p^{-m}}\big]\big(z(\tau_m - \tau_{m-1})\big).
   $$
   Thanks to \rref{1053} and the definition of $\lambda_m$, that does not exceed the value of the game with the
   payoff
   $$%\sup_{A\in\mA}\inf_{z\in A\cap \Gamma(\omega)}\Big[
   \int_{0}^{\tau_{m}-\tau_{m-1}}\lambda_m g\big(z(t)\big)\,dt+p^m
   V_*\big(z(\tau_{m}-\tau_{m-1})\big)+\varkappa.
   $$
    Then, in the last game, there exists an
     $\varkappa$-optimal strategy $A^m\in\mA$, i.e., for all $z\in A^m$, the number
$p^{m-1}V_*\big(z(0)\big)$ does not exceed
\begin{eqnarray*}
    \int_{0}^{\tau_{m}-\tau_{m-1}}\lambda_m g\big(z(t)\big)\,dt+p^m
   V_*\big(z(\tau_{m}-\tau_{m-1})\big)+3\varkappa.
\end{eqnarray*}
   Since the right-hand side of these inequalities depends only on $z|_{(0,\tau_{m}-\tau_{m-1})},$
 the strategy $A^m$ can be substituted here with any strategy $A\in\mA$
  that can be represented in the form
  $A^{m}\diamond_{\tau_m-\tau_{m-1}}A''.$
  Now, for all $A',A''\in\mA,z\in A'\diamond_{\tau_{m-1}}A^2\diamond_{\tau_m}A''$,
\begin{eqnarray}
\nonumber   p^{m-1}V_*\big(z(\tau_{m-1})\big)
    &\leq&
     \int_{\tau_{m-1}}^{\tau_m}\lambda_m g\big(z(t)\big)\,dt+p^m
   V_*\big(z(\tau_m)\big)+3\varkappa\\
   &\leqref{509}&
   e^{M/k}
   \int_{\tau_{m-1}}^{\tau_m} \varsigma'(t)g\big(z(t)\big)\,dt
   +p^m
   V_*\big(z(\tau_m)\big)+3\varkappa.
\label{555}
\end{eqnarray}

    Set $A^*\rav A^{1}\diamond_{\tau_1}A^{2}\diamond_{\tau_2}\dots \diamond_{\tau_{k-1}}A^{k^2}.$
    For every $z\in A^*$, \rref{555} holds for all $m=1,\dots,k^2.$
       Now, taking first $m=1$ and accounting for $\tau_0=0$, and subsequently applying
\rref{555} for $m=2,3\dots,k^2$, we see that
\begin{eqnarray*}
 V_*\big(z(0)\big)
    &\leqref{555}&
   e^{M/k}
   \int_{0}^{\tau_1} \varsigma'(t)g\big(z(t)\big)\,dt
   +p
   V_*\big(z(\tau_1)\big)+3\varkappa\\
   &\leqref{555}&
   e^{M/k}
   \int_{0}^{\tau_2} \varsigma'(t)g\big(z(t)\big)\,dt
   +p^2
   V_*\big(z(\tau_2)\big)+6\varkappa\\
  \dots&\leqref{555}&\dots\\
   &\leqref{555}&
   e^{M/k}
   \int_{0}^{\tau_{k^2}} \varsigma'(t)g\big(z(t)\big)\,dt
   +p^{k^2}
   V_*\big(z(\tau_{k^2})\big)+3k^2\varkappa.
 \end{eqnarray*}
for all $z\in A^*.$ Recall that $e^{M/k}<1+\epsi.$
Then, %by definitions of $\varkappa$ and $\mu'$,
we
obtain
\begin{eqnarray*}
 V_*\big(z(0)\big)
 &\leqref{M}&
    e^{M/k}\int_{0}^{\tau_{k^2}} \varsigma'(t)g\big(z(t)\big)\,dt
   +2\epsi\\
% e^{M/k}
%   \int_{0}^{\infty} \varsigma'(t)g\big(z(t)\big)\,dt
%   +p^{k^2}+3k^2\varkappa
%  \\
%
&\leq&
   e^{M/k}-1+\int_{0}^{\infty} \varsigma'(t)g\big(z(t)\big)\,dt
    +2\epsi\\
    &\leq&
    \int_{0}^{\infty} \varsigma'(t)g\big(z(t)\big)\,dt+3\epsi\qquad
    \forall z\in A^*.
 \end{eqnarray*}
 Now, $$V_*\big(z(0)\big)\leq
  \vplm[\varsigma']\big(z(0)\big)+3\epsi
   $$
   for all $z\in A^*.$ In view of \rref{1078}, we proved that
   $V_*\leq
  \vplm[\varsigma_{\mu}]+6\epsi$
  for all sufficiently small positive $\mu.$
\bo


\begin{thebibliography}{99}




\bibitem{AlvBarditrue} O.~Alvarez and  M.~Bardi,  {\em Ergodic problems in differential games,} In Advances in
    dynamic game theory, Birkh\"{a}user, Boston, (2007), pp.~131--152


\bibitem{arisconst}
 M. Arisawa   {\em Ergodic problem for the Hamilton-Jacobi-Bellman equation II}, In Annales l`Institut Henri Poincare(C) Non Linear Analysis, 15, (1998) pp.~1--24
%In Annales de l'Institut Henri Poincare (C) Non Linear Analysis (Vol. 15, No. 1, pp. 1-24). Elsevier

\bibitem{AG2000}
%\greenb{[Artstein, Gaitsgory 2000]}\\
Z. Artstein and V. Gaitsgory   {\em  The value function of singularly perturbed control systems}, Appl.
Math. Optim.,  41,3 (2000), pp.~425--445

\bibitem{BardiGame} M. Bardi,  {\em On differential games with long-time-average cost},  In Advances in dynamic games
    and their applications, Birkh\"{a}user, Boston, (2009), pp.~3--18

\bibitem{BGQ2013}
    R. Buckdahn, D. Goreac, and M. Quincampoix  {\em Existence of asymptotic values for
    nonexpansive stochastic control systems,} Appl. Math. Optim.,
    %Applied Mathematics \& Optimization,
    70,1 (2014),
    pp.~1--28


\bibitem{SIAM2012} P. Cardaliaguet, R. Laraki, and S. Sorin, {\em A Continuous Time Approach for the Asymptotic Value
    in Two-Person Zero-Sum Repeated Games,} SIAM J. Cont. Optim., 50 (2012), pp.~1573--1596

\bibitem{Cardal2010} P. Cardaliaguet  {\em Ergodicity of Hamilton-Jacobi equations with a non coercive non convex
    Hamiltonian in $\mm{R}^2/\mm{Z}^2$,} In Annales l`Institut Henri Poincare(C) Non Linear Analysis, 27,3 (2010),
    pp.~837--856.

\bibitem{QG2013} V. Gaitsgory and M. Quincampoix {\em On sets of occupational measures generated by a deterministic
    control system on an infinite time horizon,}
 Nonlinear Anal. 88, (2013), pp.~27--41.
%Nonlinear Analysis:Theory, Methods \& Applications,

%  \bibitem{Khlopin1}
%  D. V. Khlopin
%  {\em A uniform Tauberian theorem for abstract game,}
%  In: Abstracts of IV international school-seminar `Nonlinear analysis and
%extremal problems', Irkutsk, Russia, June 22Ц28, 2014. p.59


%  \bibitem{Khlopin2}   D. V.  Khlopin   {\em A uniform Tauberian theorem for conflict-controlled system,} %  In:
%Abstracts of International Conference
%  `Systems Dynamics and Control Processes' dedicated to the 90th Anniversary %
%of N.N.Krasovskii, Ekaterinburg, Russia, September 15Ц22, 2014. pp.~204--206 (in Russian)

%``''

\bibitem{KhlopinArXiv}   D. V.  Khlopin {\em  On  uniform Tauberian theorems for dynamic games,} arXiv preprint arXiv:1412.7331
    (2014)

 \bibitem{MN1981}
J. F. Mertens and A. Neyman {\em Stochastic Games,} Int. J. of Game Theory 10(1981), pp.~53--66
%Condition 3 in
%Theorem 4.1

\bibitem{1993} D. Monderer and S. Sorin {\em  Asymptotic properties in Dynamic Programming,} Int. J. of Game Theory,
    22(1993), pp.~1--11.

\bibitem{barton}
  M. Oliu-Barton, and G. Vigeral
{\em A uniform Tauberian theorem in optimal control,} In: Advances in Dynamic Games, Birkh\"{a}user, Boston, 2013,
pp.~199--215.

 \bibitem{QuinRenault}
 M. Quincampoix and J. Renault
 {\em On the existence of a limit value in some non expansive optimal
control problems,} SIAM J. Control. Optim., 49 (2011), pp.~2118Ц-2132

\bibitem{Renault2013} J. Renault {\em General limit value in Dynamic Programming,} arXiv preprint arXiv:1301.0451
    (2013)

\bibitem{Vigeral2013}
G. Vigeral  {\em A zero-sum stochastic game with compact action sets and no asymptotic value,} Dyn. Games and Appl., 3 (2013), pp.~172--186.

\bibitem{Ziliotto2013}
B. Ziliotto  {\em 
A Tauberian theorem for nonexpansive operators and applications to zero-sum stochastic games}, arXiv  preprint arXiv:1501.06525 (2015)

\end{thebibliography}
\end{document}